\documentclass[oneside,openany,a4paper]{amsart}
\usepackage[OT2,T1]{fontenc}
\DeclareSymbolFont{cyrletters}{OT2}{wncyr}{m}{n}
\DeclareMathSymbol{\Sha}{\mathalpha}{cyrletters}{"58}
\usepackage[a4paper, margin=1in]{geometry}
\usepackage{amssymb, diagrams}
\usepackage[centertags]{amsmath}
\usepackage{amsthm}
\usepackage{tikz}

\usepackage{tikz-cd}

\usepackage{amsfonts}
\usepackage{eucal}
\usepackage{mathrsfs}

\usepackage[all,cmtip]{xy}
\usepackage{graphicx, eepic}
\usepackage{verbatim}
\usepackage{lmodern}
\usepackage[utf8]{inputenc}
\usepackage[T1]{fontenc}
\usepackage{xcolor}

\usepackage{bbm}
\usepackage{pdfpages}
\usepackage{titletoc}
\usepackage{booktabs}
\usepackage[shortlabels]{enumitem}
\usepackage{mathtools}
\usepackage{thmtools}

\usepackage{epstopdf}
\usepackage{epsfig}
\usepackage{microtype}

\usepackage{hyperref}
\hypersetup{colorlinks}
\hypersetup{
    bookmarks=true,         
    unicode=false,          
    pdftoolbar=true,        
    pdfmenubar=true,        
    pdffitwindow=false,     
    pdfstartview={FitH},    
    pdftitle={My title},    
    pdfauthor={Author},     
    pdfsubject={Subject},   
    pdfcreator={Creator},   
    pdfproducer={Producer}, 
    pdfkeywords={keyword1} {key2} {key3}, 
    pdfnewwindow=true,      
    colorlinks=true,       
    linkcolor=red,          
    citecolor=magenta,        
    filecolor=violet,      
    urlcolor=blue      
}
\usepackage{multirow, url}
\usepackage{textcomp}
\DeclareMathAlphabet{\cmcal}{OMS}{cmsy}{m}{n}
\newtheoremstyle{thm}
  {3pt}
  {3pt}
  {\em}
  {0pt}
  {\bfseries}
  {}
  {5pt}
  {}
\newtheoremstyle{rem}
  {3pt}
  {3pt}
  {}
  {0pt}
  {\bfseries}
  {.}
  {5pt}
  {}

\newtheorem{thm}{Theorem}[section]

\newtheorem{prop}[thm]{Proposition}

\newtheorem{conj}[thm]{Conjecture}

\theoremstyle{definition}

\theoremstyle{rem}
\newtheorem{rem}[thm]{{Remark}}

\numberwithin{equation}{section} \numberwithin{table}{section}
	
\newtheorem*{thm*}{Theorem}
\newtheorem*{rem*}{Remark}
\newtheorem*{rems*}{Remarks}
\newtheorem*{exam*}{Example}
\newtheorem*{exams*}{Examples}

\makeatletter
\newcommand{\neutralize}[1]{\expandafter\let\csname c@#1\endcsname\count@}
\makeatother

\newcommand{\Z}{{\mathbb{Z}}}

\newcommand{\fg}{{\mathfrak{g}}}

\newcommand{\cA}{{\cmcal{A}}}
\newcommand{\cB}{{\cmcal{B}}}

\newcommand{\cM}{{\cmcal{M}}}

\newcommand{\cO}{{\cmcal{O}}}

\def\n2Z{\frac{1}{n^2}\Z/\Z}

\def\d{\delta}

\newcommand{\Gm}{{\mathbb{G}}_{m}}

\newcommand{\arinj}{\ar@{^(->}}
\newcommand{\arsurj}{\ar@{->>}}
\newcommand{\arsub}{\ar@{}[r]|-*[@]{\subset}}
\newcommand{\arsup}{\ar@{}[r]|-*[@]{\supset}}
\newcommand{\arcap}{\ar@{}[d]|-*[@]{\subset}}
\newcommand{\arcup}{\ar@{}[u]|-*[@]{\subset}}
\newcommand{\arin}{\ar@{}[u]|-*[@]{\in}}

\renewcommand{\Im}{\operatorname{Im}}

\newcommand{\Spec}{{\operatorname{Spec}}}

\renewcommand{~}{\hspace{0.4mm}}

\newcommand{\bd}{\boldsymbol}

\mathchardef\hyp="2D

\def\<{\langle }

\def\>{\rangle}

\def\nZ{\frac{1}{n}\mathbb{Z}/\mathbb{Z}}
\def\d{\delta}
\def\ms{\medskip}

\def\nZ{\frac{1}{n}\Z/\Z}

\def\UX{\cO_X^{\times}}
\newcommand{\ZZ}{{\mathbb Z}}

\DeclareMathOperator{\divis}{div}
\DeclareMathOperator{\Div}{Div}
\DeclareMathOperator{\Cl}{Cl}

\def\Om{\Omega}
\def\bd{\begin{diagram}}
\def\ed{\end{diagram}}

\begin{document}                                                                          

\title{A note on abelian arithmetic BF-theory}
\author{Magnus Carlson}
\address{M.C:  Department of Mathematics, Hebrew University}
                                                     
\author{Minhyong Kim}
\address{M.K: Mathematical Institute, University of Oxford, and the Korea Institute for Advanced Study}

\maketitle

\begin{abstract}
We compute some arithmetic path integrals for $BF$-theory over the ring of integers of a totally imaginary field, which evaluate to natural arithmetic invariants associated to $\Gm$ and abelian varieties.
\end{abstract}


\section{Towards arithmetic BF theory}
$BF$-theory is a rare example of a topological field theory that can be defined in any dimension \cite{Cattaneo}. Let $M$ be an oriented $n$-manifold, $G$ a compact Lie group, and $P$ a principal $G$-bundle on $M$. Assume given also a finite-dimensional representation $\rho$ of the Lie algebra $\fg$ of $G$. The $BF$ functional depends on two fields, a connection $A$ on $P$, and  an $(n-2)$-form $B\in \Om^{n-2}(M,adP)$ with values in the adjoint bundle of $P$. Then
$$BF(B,A):=\int_MTr[\rho(B)\wedge \rho(F_A)],$$
where $F_A$ is the curvature form of $A$. (There seem to be some other possibilities for the invariant function on $\fg$ that goes into the integral.) A number of properties of $BF$-theories make them easier to deal with than Chern-Simons theories. Primary among them for our purposes is that the two variables allow us to extend the definition to non-orientable manifolds, simply by letting the field $B$ vary over $\Om^{n-2}(M,adP\otimes \omega_M)$, where $\omega_M$ is the orientation bundle of $M$. 

In some earlier papers, a preliminary attempt to define and compute arithmetic analogues of Chern-Simons functions was made \cite{KimLinking, KimCSI, KimCSII, AhlquistCarlson, Bleher}. Also, moduli spaces of `arithmetic gauge fields' have been applied to Diophantine geometry \cite{Balakrishnan, KimGauge}. One of the obstructions  to developing a full-fledged  arithmetic topological field theory based on Chern-Simons theory is that  natural arithmetic dualities involve a sheaf $\mu_n$  or $\hat{\Z}(1)$, which are not trivialisable in general.  That is, arithmetic schemes are not orientable. There are various ways to circumvent this problem, all of which introduce some difficulties for working out interesting examples. The purpose of this paper is to suggest, by way of two brief computations, that BF-theory provides a simpler way to link number theory to topological field theory. 

Let $X=\Spec(\cO_F)$, where $F$ is a totally imaginary number field. Then
\begin{prop}
$$\sum_{(a,b)\in H^1(X, \Z/n)\times H^1(X, \mu_{n})} \exp(2\pi i BF(a,b))=|n\Cl_F[n^2]|   \cdot | \UX/(\UX)^{n}| \cdot |\Cl_F/n|.$$ 
\end{prop}
Thus, for $n$ large enough, the `finite path integral' will capture exactly $$| \UX/(\UX)^{n}| \cdot |\Cl_F|,$$a quantity of the form (regulator $\times$ class number). The precise notation and definitions that go into this proposition as well as the next one will be explained in subsequent sections.

Now let $A$ and $B$ be dual abelian varieties over $F$ with semi-stable reduction at all places. Let $n$ be an integer coprime to all the Tamagawa numbers of $A$ and $B$ as well as to the places of bad reduction for them. Assume that the Tate-Shafarevich groups $\Sha(A)$ and $\Sha(B)$ are finite  and $\Sha(B)[n]=\Sha(B)[n^2]$. Then
\begin{prop}
$$\sum_{(a,b)\in H^1(X,\cA[n])\times H^1(X,\cB[n])} \exp(2\pi i BF(a,b))=|A(F)/n| \cdot |B(F)/n| \cdot |\Sha(A)[n]|$$
\end{prop}
Here, $\cA$ and $\cB$ are the Neron models of $A$ and $B$ respectively. Note that the expression is symmetric in $A$ and $B$ because $|\Sha(A)[n]|=|\Sha(B)[n]|$ via the Cassel-Tate pairing.

We view these formulae as some (weak) evidence for the suggestion made in \cite{KimCSI} that an arithmetic topological functional on moduli spaces of Galois representations will have something to do with $L$-functions. In any case,  it is rather striking to find expressions for the orders of class groups and Tate-Shafarevich groups as exponential sums, which perhaps have not appeared heretofore  in the literature. Equally  notable is that a path integral for the BF functional for three manifolds leads to the Alexander polynomial of a knot in the physics literature \cite{Cattaneo}, which is well-known to be an analogue of the $p$-adic $L$-function \cite{Morishita}.

In forthcoming work, we will develop  $BF$-theory for arithmetic schemes with boundary, that is $\Spec(\cO_F[1/S])$ for a finite set $S$ of places, as well as a $p$-adic theory. Also interesting would be to develop arithmetic $BF$-theory for arithmetic schemes in higher dimension by way of the duality theory of \cite{Geisser}. However, in the  present announcement, the primary goal is to illustrate  with a minimum of clutter the relationship between a path integral in the sense of physicists and important arithmetic invariants.

\section{Some finite path integrals for $\mathbb{G}_m$: Proof of Proposition 1.1.}

As before, let $X=\Spec(\cO_F)$ for $F$ a totally imaginary number field and let $\mu_{n}$ and $\Z/n$ be the usual finite flat group schemes viewed as sheaves in the flat topology. We have the Bockstein map
$$\d: H^1(X, \mu_{n})\rTo H^2(X,\mu_{n})$$
coming from the exact sequence
$$1\rTo\mu_{n}\rTo \mu_{n^2}\rTo \mu_{n}\rTo 1$$
and the invariant isomorphism \cite{MazurNotes}
$$\int :H^3(X, \mu_{n})\rTo \nZ.$$
Define the BF-functional on $$\cM:=H^1(X, \Z/n)\times H^1(X, \mu_{n})$$ by
$$BF(a,b)=\int(a \cup \d b).$$
(The class $\d b$ is the analogue of the curvature $F$.)
\begin{rem}
One can also define a BF-functional on $\cM$ as $BF'(a,b) = \int(\d a \cup b),$ where $\d:H^1(X,\ZZ/n) \rightarrow H^2(X,\ZZ/n^2)$ is the Bockstein map coming from the exact sequence $$1 \rTo \ZZ/n \rTo \ZZ/n^2 \rTo \ZZ/n \rTo 1.$$ However, $BF' = BF,$ since we have an equality $\d a \cup b = a \cup \d b$ by Lemma 2.1 and the proof of Lemma 2.2 in \cite{KimLinking}.
\end{rem}
We will now calculate the path integral
$$\sum_{(a,b)\in \cM} \exp(2\pi i BF(a,b)).$$
First, let us calculate the groups $H^i(X,\mu_n).$  We define $\Div F$ to be the group of fractional ideals of $F$ and for $x \in F^*,$ we let $\divis(x)$ be the associated principal ideal. We claim that $$H^i(X,\mu_n) = \begin{cases} \mu_n(F) & \text{for } i=0  \\ Z_1/B_1 & \text{for } i = 1 \\ \Cl_F/n & \text{for } i=2 \\ \ZZ/n & \text{for } i=3 \\ 0 & \text{for } i >3 \end{cases}$$ where $$Z_1 = \{ (x,I) \in F^* \oplus \Div F | \ nI = -\divis(x)\}$$ and $B_1 = \{(x^n,-\divis(x))\in F^* \oplus \Div F | \ x \in F^*\}.$ To see this, note that $\mu_n,$ is quasi-isomorphic to the complex $$\mathcal{C}_\bullet = (\mathbb{G}_{m,X} \xrightarrow{\cdot n} \mathbb{G}_{m,X})$$ with non-zero terms in degree $0$ and $1.$ Since $\mathbb{G}_{m,X}$ is representable by a smooth group scheme, flat cohomology groups with coefficients in $\mathbb{G}_{m,X}$ coincide with the corresponding étale cohomology groups \cite[III, Theorem 3.9]{MilneEtale}. Thus, to compute flat cohomology with coefficients in $\mu_n,$ it suffices to compute the hypercohomology of $\mathcal{C}_\bullet,$ viewed as a complex of sheaves on the étale site. To do this, we take the well-known resolution $j_* \mathbb{G}_{m,F} \xrightarrow{\divis} \Div X$ of $\mathbb{G}_{m,X}.$ Here, $$j:\Spec F \rightarrow X$$ is the inclusion of the generic point and $$\Div X = \underset{x \in X^0}{\oplus} \ZZ_{/x},$$ where $\ZZ_{/x}$ is the skyscraper sheaf at the closed point $x.$ Resolving the complex $\mathcal{C}_\bullet$ levelwise with this resolution and taking the total complex of the resulting double complex, we get the complex 
$$j_* \mathbb{G}_{m,F} \xrightarrow{(\divis,n^{-1})} \Div X \oplus j_* \mathbb{G}_{m,F} \xrightarrow{n+\divis} \Div X.$$ Now, as in Lemma 4.2 and Corollary 4.3 of \cite{CarlsonSchlank} one uses the local-to-global spectral sequence to conclude that the complex above computes the flat cohomology of $\mu_n$ for $i= 0,1,2;$ to conclude the calculation for  $i \geq 3$ one applies flat duality as in III, Corollary 3.2 of \cite{MilneADT}.
We now proceed to analyze the Bockstein map $$\d: H^1(X,\mu_n) \rTo H^2(X,\mu_n).$$ Using the complex given earlier, the Bockstein was computed in Lemma 4.1 of \cite{AhlquistCarlson}, and it is the composite of two maps: the first is the surjective map which takes $(x,I) \in H^1(X,\mu_n) \cong Z_1/B_1$ to $I \in \Cl_F[n]$ and the second is the reduction map $\Cl_F[n] \rightarrow \Cl_F/n.$  By noting that the kernel of the first map is $\UX/(\UX)^n$ and that the kernel of the second map is $n \Cl_F[n^2],$ we see that our sum becomes  $$|\UX/(\UX)^n| \cdot |n \Cl_F[n^2]| \sum_{(a,b) \in H^1(X,\ZZ/n\ZZ) \times \Cl_F[n]/n\Cl_F[n^2]} \exp(2 \pi i a \cup \bar{b})$$ where $\bar{b} \in \Cl_F/n$ is the reduction of $b.$ But for $b$ non-trivial, it is clear that this sum is zero, giving us
$$\sum_{(a,b)\in \cM} \exp(2\pi i BF(a,b))=|n\Cl_F[n^2]|   \cdot | \UX/(\UX)^{n}| \cdot |\Cl_F/n|.$$ In particular, if $\Cl_F[n] = \Cl_F,$ we see that the sum evaluates to $|\UX/(\UX)^{n}| \cdot |\Cl_F|.$ If one calculates the path integral $\sum_{(a,b) \in \cM} \exp(2\pi i BF'(a,b))$ where $BF'$ is as in the above remark, one sees, since $BF'$ coincides with $BF,$ that $\sum_{(a,b) \in \cM} \exp(2 \pi i BF(a,b))$ equals $|\Cl_F[n]| \cdot |\UX/(\UX)^n|$ times the number of isomorphism classes of unramified $\ZZ/n$-étale algebras which can be embedded into an $\ZZ/n^2$-étale algebra.
\ms
\section{Some finite path integrals for abelian varieties: Proof of Proposition 1.2}
 Let $\cA$ be the Néron model of an abelian variety $A$  over $F$ and let $\cB$ be the Néron model of the dual $B$. Assume that both $\cA$ and $\cB$ have semi-stable reduction and let $n$ be a sufficiently large positive integer. More precisely, suppose both $A$ and $B$ have good reduction at all places dividing $n$. Further assume that $n$ is coprime to $|\Phi_A |\cdot |\Phi_B|$, where $\Phi_A$ and $\Phi_B$ are the groups of connected components of $\cA$ and $\cB$. According to \cite[Theorem 1.1(a)(ii)]{Ces}, we have isomorphisms 
$$H^1(X, \cA[n])\cong Sel(F,A[n]),$$
$$H^1(X, \cB[n])\cong Sel(F,B[n])$$
where the left hand side is flat cohomology (in the fppf site).
We claim that our assumption that $\cB$ has semi-stable reduction, together with our assumptions on $n$, implies that multiplication by $n$ on $\cB$ is an epimorphism in the category of fppf sheaves. To see this, note that by \cite[Lemma B.4.]{Ces}, $\cB^0 \xrightarrow{\cdot n} \cB^0$ is faithfully flat, thus an epimorphism of fppf sheaves. Further, we have a commutative diagram
$$\begin{tikzcd} 
0 \arrow[r] & \cB^0 \arrow[r] \arrow[d,"n"] & \cB \arrow[r] \arrow[d,"n"] & \Phi_B \arrow[r] \arrow[d,"n"] & 0 \\
0 \arrow[r] & \cB^0 \arrow[r] & \cB  \arrow[r] & \Phi_B \arrow[r] & 0  
\end{tikzcd}$$
and since $n$ is prime to $|\Phi_B|,$ multiplication by $n$ on $\Phi_B$ is actually an isomorphism. Thus, by the snake lemma, multiplication by $n$ on $\cB$ is an epimorphism. By \cite[Corollary 3.4]{MilneADT} there is a perfect pairing 
$$\cup : H^1(X,\cA[n]) \times H^2(X,\cB[n]) \rightarrow H^3(X, \mathbb{G}_m) \cong \mathbb{Q}/\mathbb{Z}.$$ 
Using the exact sequence
$$0\rTo \cB[n]\rTo \cB[n^2]\rTo^{n} \cB[n]\rTo 0,$$
we get the Bockstein map
$$\d: H^1(X,\cB[n])\rTo H^2(X,\cB[n]).$$ We now define the BF-functional
$$BF: H^1(X,\cA[n])\times H^1(X,\cB[n])\rTo H^3(X,\mathbb{G}_m)[n] \cong H^3(X,\mu_{n}) \cong \dfrac{1}{n}\ZZ/\ZZ$$
$$BF(a,b) = \int (a\cup \d b),$$ where $\int$ is the isomorphism $H^3(X,\mu_{n}) \cong \dfrac{1}{n}\ZZ/\ZZ,$ as in the previous section. Our goal is now to calculate the path integral $$\sum_{(a,b)\in H^1(X,\cA[n])\times H^1(X,\cB[n])} \exp(2\pi i BF(a,b)).$$

To this end, note that, by \cite[Proposition 4.5 (c)-(d)]{Ces} and our assumptions on $n,$ we have
$$H^1(X,\cA)[n]\cong \Sha(A)[n], \ \ H^1(X,\cB)[n]\cong \Sha(B)[n].$$ 
Assuming the finiteness of the Tate-Shafarevich group, we can take $n$ large enough so that $$H^1(X,\cB)[n^2]=H^1(X,\cB)[n].$$

We now identify the kernel of $\d.$ By using the diagram
$$\bd 
0&\rTo & \cB/n^2\cB&\rTo& H^1(X,\cB[n^2])&\rTo &H^1(X,\cB)[n^2]&\rTo 0\\
 && \dTo& & \dTo^{n}& &\dTo^{n}& \\
 0&\rTo & \cB/n\cB&\rTo& H^1(X,\cB[n])&\rTo &H^1(X,\cB)[n]&\rTo  0\ed$$
and the fact that the rightmost map is zero, we see that $Ker(\d)=Im(\cB/n\cB).$ Thus, we get an induced injection
$$\bar{\d}:H^1(X,\cB)[n]\rInto H^2(X,\cB[n]).$$ 
By an argument identical to the case of $\mu_{n}$ above, since $$|H^1(X,\cA[n])| =  |A(F)/n| \cdot |\Sha(A)[n]|,$$ we get
$$\sum_{(a,b)\in H^1(X,\cA[n])\times H^1(X,\cB[n])} \exp(2\pi i BF(a,b))=|B(F)/n| \sum_{(a,c)\in H^1(X,\cA[n])\times Im(\d)}\exp(2\pi i a\cup c)$$

$$=|B(F)/n| \sum_{(a,c)\in H^1(X,\cA[n])\times \Im(\d)}\exp(2\pi i a\cup c)$$

$$=|A(F)/n| \cdot |B(F)/n| \cdot |\Sha(A)[n]|.$$

\section*{Acknowledgements}
M.K. is grateful to Kevin Costello and Edward Witten for urging him to look at BF-theory, and to Nima Arkani-Hamed for an invitation to speak at the IAS high-energy theory seminar in the course of which the initial ideas for an arithmetic BF-theory came to mind. He is also grateful to Philip Candelas, Xenia de la Ossa, Tudor Dimofte, Rajesh Gopakumar, Sergei Gukov,  Jeff Harvey, Theo Johnson-Freyd, Albrecht Klemm, Si Li, Tony Pantev, Ingmar Saberi, Johannes Walcher, Katrin Wendland, and Philsang Yoo for numerous illuminating conversations on quantum field theory. He was supported in part by the EPSRC programme grant EP/M024830/1, `Symmetries and Correspondences'. 

M.C. would like to thank the  Wallenberg foundation for their support. He is also grateful to Merton College for a visiting research scholarship in April of 2019, thanks to which it was possible to begin the research presented in this paper.


\begin{thebibliography}{10}
\bibitem{AhlquistCarlson} {Ahlquist, E., Carlson, M.,The cohomology ring of the ring of integers of a number field. arxiv:1612.01766} 
\bibitem{Balakrishnan}{Balakrishnan, J.,  Dan-Cohen, Kim, M., I., Wewers, S.) A non-abelian conjecture of Tate-Shafarevich type for hyperbolic curves.
Mathematische Annalen, October 2018, Volume 372, Issue 1-2, pp 369--428.}
\bibitem{Bleher}{Bleher, F.,  Chinburg, T.,  Greenberg, R.,  Kakde, M.,  Pappas,  G., Martin J. Taylor, M., Unramified arithmetic Chern-Simons invariants. arXiv:1705.07110v1.

}

\bibitem{CarlsonSchlank} {Carlson, M., Schlank, T., The unramified inverse Galois problem and cohomology rings of totally imaginary number fields. arxiv:1612.01766} 
\bibitem{Cattaneo}{ Cattaneo, A.S.,  Cotta-Ramusino, P., Fr\"ohlich, J., Martellini, M., Topological BF theories in 3 and 4 dimensions.
Journal of Mathematical Physics 36, 6137 (1995)}
\bibitem{Ces} {Cesnavicius, K., Selmer groups as flat cohomology groups. arXiv:1301.4724}
\bibitem{KimLinking} {Chung, H., Kim, D., Kim, M., Pappas, G., Park, J., Yoo, H., Abelian Arithmetic Chern–Simons Theory and Arithmetic Linking Numbers, arxiv:1706.03336}
\bibitem{Geisser}{Geisser, T., Schmidt, A., Poitou-Tate duality for arithmetic schemes. Compositio Math. 154 (2018) 2020-2044}
esearch Notices, Published electronically, November, 2017
\bibitem{KimCSI}{Kim, M., Arithmetic Chern-Simons Theory I. arXiv:1510.05818}
\bibitem{KimCSII}Chung, H.,  Kim, D., Kim, M., Park, J., Yoo, H., Arithmetic Chern-Simons Theory II.
(with ) To be published in p-adic Hodge Theory, Simons Symposia, Springer-Verlag.
\bibitem{KimGauge}{Arithmetic Gauge Theory: A Brief Introduction
Modern Physics Letters A, Volume 33, Issue 29 (2018).}
\bibitem{MazurNotes} {Mazur, B., Notes on étale cohomology of number fields, Ann.Sci. ÈcoleNorm.Sup. (4) 6, 521-552,1973}
\bibitem{MilneADT} {Milne, J.S, Arithmetic duality theorems, BookSurge, LLC, Charleston, SC, Second edition, 2006}
\bibitem{MilneEtale} {Milne, J.S, Étale cohomology, volume 33 of Princeton Mathematical Series. Princeton University Press, Princeton, N.J., 1980.}
\bibitem{Morishita}{Morishita, M, {\em Knots and Primes}.  Springer-Verlag (2012).}





\end{thebibliography}
\end{document}